\documentclass[11pt]{article}   
\usepackage{amssymb,amscd,latexsym}   
\usepackage{amsmath}
\usepackage{amsthm}
\newcommand{\rar}{\rightarrow}

\newcommand{\llar}{-\kern-5pt-\kern-5pt\longrightarrow}

\newtheorem{Theorem}{Theorem}[section]

\newtheorem{Proposition}[Theorem]{Proposition}

\def\sqr#1#2{{\vcenter{\hrule height.#2pt
        \hbox{\vrule width.#2pt height#1pt \kern#1pt
            \vrule width.#2pt}
        \hrule height.#2pt}}}
\def\phi{\varphi}
\def\demo{\noindent{\bf Proof. }}
\def\square{\mathchoice\sqr64\sqr64\sqr{4}3\sqr{3}3}
\def\qed{\hspace*{\fill} $\square$}

\def\XX{{\bf X}}
\def\YY{{\bf Y}}

\def\ZZ{{\bf Z}}

\def\hht{{\rm ht}\,}

\def\ass{{\rm Ass}\,}

\def\restr{{\kern-1pt\restriction\kern-1pt}}

\def\pp{{\mathbb P}}

 2
 3
 4
 5
 4
 3
 2
 1
 1
 2
 3
 1
 0
 1
 2
 3
 4
 5

\begin{document}
\begin{center}
{\Large{\bf\sc On catalecticant perfect ideals of codimension $2$}}
\footnotetext{2010 AMS {\it Mathematics Subject
Classification}: 13A30, 13C14, 13C40, 13D02, 13E15, 13H10, 14E05, 14E07.}

\vspace{0.3in}

{\large\sc Zaqueu Ramos}\footnote{Parts of this work were done while this author held a
Doctoral Fellowship (CAPES, Brazil).} \quad\quad
 {\large\sc Aron  Simis}\footnote{Partially
supported by a CNPq grant and a PVNS Fellowship (CAPES).}

\end{center}


\bigskip

\begin{abstract}

One deals with catalectic codimension two perfect ideals
and certain degenerations thereof, with a view towards the nature of their symbolic powers.
In the spirit of \cite{zaron} one considers linearly presented such ideals,
only now in the situation  where the number of variables
is sufficiently larger than the size of the matrix, yet still stays within reasonable bounds.
\end{abstract}

\section*{Introduction}

The idea in this work is to consider a codimension $2$ perfect ideal $I$ in a polynomial ring
$R=k[X_1,\ldots,X_n]$ ($k$ a field) whose $m\times (m-1)$ presentation matrix
is a generic catalecticant or a degeneration of this in a sense to be explained.
We examine three of these constructs: catalecticants with leap, sub-Hankel and semi-Hankel matrices.
Depending on the values of $n$ and $m$ in the sub-Hankel case the ring $R/I$ may fail to be a normal domain
(this is to be compared to with \cite[Proposition 2.3]{zaron}, where normality always holds if $n\geq 4$) .
As it turns out, a special instance of the semi-Hankel case leads to a Cremona transformation.
Such a Cremona transformation does not seem to have been observed before in this systematization.

As a rule, the sort of results obtained for all the above constructs have to do with the properties of the
ideals of minors, lower order minors as well.
Since mostly $m<<n$, there is room enough to ask whether $I$ is of linear type or at least locally so
on the punctured spectrum, typically via variations of the so-named $(G_s)$ condition.
The other property of interest is normally tosionfreeness.
Mostly, $I$ will be normally torsionfree, but there are some exceptions.
It becomes clear that whether this property holds has to do with the nature of the degeneration
of the catalecticant (or the generic matrix).

Since basically every section deals with one of the constructs and there is essentially one theorem in each
case, there is hardly any point in enlarging this foreword with further details.
As a way of compensation, we will expand on the terminology and tools used throughout.

For the proofs we have drawn quite a bit on the results of \cite{Eisenbud2}, which take a central
role in this work.

\section{Terminology}

We will assume throughout that $R$ is a standard graded polynomial ring over an infinite field $k$.
Given an ideal $I\subset R$ and an integer $r\geq 1$, the $r$th {\em symbolic power} $I^{(r)}$ of $I$
is the contraction of $U^{-1}I^{r}$ under the natural homomorphism $R\rar U^{-1}R$
of fractions, where $U$ is the complementary set of the union of the associated primes of $R/I$.
In this work $I$ will be a codimension $2$ perfect ideal, hence $R/I$ is Cohen--Macaulay and so $I$ is a
pure (unmixed) ideal.
In this setup then $I^{(r)}$ is precisely the intersection of the primary components
of the ordinary power $I^{r}$ relative to the associated primes of $R/I$, i.e., the unmixed part of $I^{r}$.

One says that $I$ is {\em normally torsionfree} provided $I^{(r)}=I^r$ for every $r\geq 1$.

We observe that, like the ordinary powers, the symbolic powers constitute a decreasing multiplicative filtration,
so one can consider the corresponding {\em symbolic Rees algebra} $\mathcal{R}_{R}^{(I)}=\bigoplus_{r\geq 0}I^{(r)}t^{r}
\subset R[t]$.
However, unlike the ordinary Rees algebra, this algebra may not be finitely generated over $R$.

\medskip

Next is a review of the notion of the inversion factor associated to a Cremona map.
This idea has been largely addressed in \cite{zaron}.

Let $k$ denote an arbitrary infinite field.
A rational map $\mathfrak{G}:\pp^{n-1}\dasharrow \pp^{m-1}$ is defined by $m$ forms $\mathbf{g}=\{g_1,\ldots, g_m\}
\subset R:=k[\XX]=k[X_1,\ldots,X_n]$ of the same degree $d\geq 1$, not all null.
We naturally assume throughout that $n\geq 2$.
We often write $\mathfrak{G}=(g_1:\cdots :g_m)$ to underscore the projective setup
and assume that $\gcd\{g_1,\cdots ,g_m\}=1$ (in the geometric terminology, the linear system defining $\mathfrak{G}$ ``has no fixed part''),
in which case we call $d$ the {\em degree} of $\mathfrak{G}$.

We say that $\mathfrak{G}$ is a {\em Cremona map} if $m=n$ and $\mathfrak{G}$  is a birational map of $\pp^{n-1}$.
This means that there  is a rational map
 $\pp^{n-1}\dasharrow \pp^{n-1}$ with defining coordinate forms $\mathbf{f}=\{f_1,\ldots, f_n\}
\subset k[\YY]$ satisfying the
relations
\begin{equation}\label{birational_rule}\nonumber
({f}_1(\mathbf{g}):\cdots :{f}_n(\mathbf{g}))=(X_1:\cdots :X_n), \;
({g}_1(\mathbf{f}):\cdots :{g}_n(\mathbf{f}))= (Y_1:\cdots :Y_n)
\end{equation}
The first of the above structural congruences
\begin{equation}\label{congruence}
(f_1(g_1,\ldots,g_m),\ldots, f_n(g_1,\ldots,g_n))\equiv (X_1,\ldots,X_n)
\end{equation}
involving the inverse map gives a uniquely defined
form $D\in R$ up to a nonzero scalar in $k$, such that $f_i(g_1,\ldots,g_m)=X_iD$, for every $i=1,\ldots,n$.

We call $D$ the {\em source inversion factor} of $\mathfrak{G}$

\medskip

The classical theory of plane Cremona maps in characteristic zero relates the Jacobian of a homaloidal net with the
principal curves of the corresponding Cremona map.
In this connection the following general result  has been proved in \cite{zaron}.

\begin{Proposition}\label{inversion_factor_is_determinant}{\rm (char$(k)=0$)}
Let $R=k[X_{1},\ldots, X_{n}]$ be a polynomial ring over a field $k$ of characteristic zero, with its standard grading
and let $\mathcal{L}=(\ell_{ij})$ be an $n\times (n-1)$ matrix whose entries are linear forms in $R$.
For every $i=1,\ldots, n$ write $\Delta_{i}$ for the signed $(n-1)$-minor of $\mathcal{L}$ obtained by omitting
the $i$-th row and let $\Theta=\Theta(\mathbf{\Delta})$ denote the Jacobian matrix of $\mathbf{\Delta}:=
\{\Delta_{1},\ldots,\Delta_{n}\}$.

If the ideal $I_{n-1}(\mathcal{L}):=(\mathbf{\Delta})\subset R$ is of linear type
then the rational map $\pp^{n-1}\dasharrow \pp^{n-1}$ defined by $\mathbf{\Delta}$
is a Cremona map and the associated source inversion factor is $\frac{1}{n-1}\det (\Theta)$.
\end{Proposition}

\section{Main results}

Following common usage, one denotes by $I_t(\Psi)\subset R$ the ideal generated by the $t\times t$ minors of a matrix $\Psi$.

\subsection{Generic catalectic matrices with leap}

The basic structure in this part is an $m\times (m-1)$ $r$-{\em leap catalecticant} (or, shortly,
$r$-{\em catalecticant}) in $R=k[X_1,\ldots,X_n]$, where $1\leq r\leq m-1$
and $n=(m-1)(r+1)$:
$$\mathcal{C}=
\begin{pmatrix}
X_1 & X_2 & X_3 & \hdots & X_{m-1}\\
X_{r+1} & X_{r+2} & X_{r+3} & \hdots & X_{m+r-1}\\
X_{2r+1} & X_{2r+2} & X_{2r+3} & \hdots & X_{m+2r-1}\\
\vdots & \vdots & \vdots & \ddots & \vdots\\
X_{(m-1)r+1} & X_{(m-1)r+2} & X_{(m-1)r+3} & \hdots & X_{(m-1)r+(m-1)}
\end{pmatrix}
$$
The extreme values $r=1$ and $r=m-1$ yield, respectively, the ordinary Hankel matrix and the generic matrix.

A crucial property of $r$-catalectic matrices is essentially contained in \cite{Eisenbud2}:

\begin{Proposition}\label{cat_1-generic}
An $r$-catalectic matrix of arbitrary size $v\times w \, (v\geq w)$ is $1$-generic.
\end{Proposition}
\demo One can embed a catalectic matrix $\mathcal{C}$ such as the one above, of arbitrary size $v\times w$,
as a submatrix of a $v'\times w$ Hankel matrix for a suitable $v'\geq v$.
Namely
$$\mathcal{C}=
\begin{pmatrix}
X_1 & X_2 & X_3 & \hdots & X_{w}\\
X_{r+1} & X_{r+2} & X_{r+3} & \hdots & X_{r+w}\\
X_{2r+1} & X_{2r+2} & X_{2r+3} & \hdots & X_{2r+w}\\
\vdots & \vdots & \vdots & \ddots & \vdots\\
X_{(v-1)r+1} & X_{(v-1)r+2} & X_{(v-1)r+3} & \hdots & X_{(v-1)r+w}
\end{pmatrix}
$$
can be augmented to a matrix $\mathcal{C}'$ with the same number of columns $w$ by adding the
``missing'' Hankel blocks
$$\mathcal{C}'=
\begin{pmatrix}
X_1 & X_2 & X_3 & \hdots & X_{w}\\
X_{r+1} & X_{r+2} & X_{r+3} & \hdots & X_{r+w}\\
X_{2r+1} & X_{2r+2} & X_{2r+3} & \hdots & X_{2r+w}\\
\vdots & \vdots & \vdots & \ddots & \vdots\\
X_{(v-1)r+1} & X_{(v-1)r+2} & X_{(v-1)r+3} & \hdots & X_{(v-1)r+w}\\[3pt]
\hline\\[-3pt]
X_2 & X_3 & X_4 &\hdots  & X_{w+1}\\
X_3 & X_4 & X_5 &\hdots  & X_{w+2}\\
\vdots & \vdots & \vdots & \ddots & \vdots\\
X_{r} & X_{r+1} & X_{r+2} & \hdots & X_{r+w-1}\\[3pt]
\hline\\[-3pt]
X_{r+2} & X_{r+3} & X_{r+4} &\hdots  & X_{r+w+1}\\
X_{r+3} & X_{r+4} & X_5 &\hdots  & X_{r+w+2}\\
\vdots & \vdots & \vdots & \ddots & \vdots\\
X_{2r} & X_{2r+1} & X_{2r+2} & \hdots & X_{2r+w-1}\\[3pt]
\hline\\[-3pt]
\vdots & \vdots & \vdots & \ddots & \vdots
\end{pmatrix}
$$
Now, given $A\in {\rm GL}_k(v)$ we consider the block matrix
$$A'=\left(
\begin{array}{c|c}
A & \mathbf{0}\\
\mathbf{0} & \mathbf{1}_{v'-v}
\end{array}
\right),
$$
where $\mathbf{1}_{v'-v}$ denotes the identity matrix os size $v'-v$.
Then $A'\in{\rm GL}_k(v')$.
This shows that if performing row operations on $\mathcal{C}'$ produces no zero entry, then the same holds for
$\mathcal{C}.$
The same argument for column operations is trivial since we have not changed the column size of the
original matrix.

Next, on $\mathcal{C}'$ we move up the additional Hankel blocks in such a way so as to have the resulting matrix become
a Hankel matrix of size $v'\times w$.
Note that this operation consists of iterated row permutations, so $\mathcal{C}'$ is $1$-generic
if and only the resulting matrix is.
Since a Hankel matrix is $1$-generic (\cite[Proposition 4.2]{Eisenbud2}), we are done.
\qed

\medskip

Using this, we prove our result of this subsection:

\begin{Theorem}\label{Perfect_catalecticant}
Let $I\subset R=k[X_1,\ldots,X_n]$ stand for the ideal of $(m-1)$-minors of an $m\times (m-1)$ $r$-leap catalectic
matrix $\mathcal{C}$ as above, with $1\leq r\leq m-1$.
Then
\begin{enumerate}
\item[{\rm (a)}] $\hht(I_t(\mathcal{C}))\geq m-t+2$ for $1\leq t\leq m-2$ and $\hht(I)=2$.
\item[{\rm (b)}] $R/I$ is a Cohen--Macaulay normal domain.
\item[{\rm (c)}] $I$ is an ideal of linear type.
\item[{\rm (d)}] $I$ is normally torsionfree.
\end{enumerate}
\end{Theorem}
\demo
(a) The result is clear for $t=1$, hence assume that $2\leq t\leq m-1$.
For $t$ in this interval, consider the submatrix $[\mathcal{C}]_t$ of $\mathcal{C}$ formed by its first $t$ columns.
By Proposition~\ref{cat_1-generic}, $[\mathcal{C}]_t$ is $1$-generic, hence
its ideal $I_t([\mathcal{C}]_t)$ of $t$-minors (maximal minors) of $[\mathcal{C}]_t$ is prime and satisfies
$\hht (I_t([\mathcal{C}]_t))\geq m-1-t+2=m-t+1$ (cf. \cite[Theorem 2.1]{Eisenbud2}).

To conclude, it suffices to show that the inclusion $I_t([\mathcal{C}]_t)\subset I_t(\mathcal{C})$ is proper for $t\leq m-2$.
For this, let $\Delta$ stand for the lower rightmost $t$-minor of $\mathcal{C}$.
Since $\Delta$ has a term involving effectively the last variable, while  no $t$-minor minor of $[\mathcal{C}]_t$
(with $t\leq m-2$) has such a term, and since all minors in consideration live in the same degree, we
clearly have $\Delta\notin I_t([\mathcal{C}]_t)$.

\smallskip

(b) Cohen--Macaulayness is obvious.
Since $\mathcal{C}$ is $1$-generic
then every prime $Q\subset R/I$ such that $(R/I)_Q$ is not regular must contain the ideal
$I_{m-2}(\mathcal{C})/I$ (\cite[Corollary 3.3]{Eisenbud2}).
Then $\hht(Q)\geq\hht(I_{m-2}(\mathcal{C})/I)\geq 4-2= 2$, by (a).
Therefore, $R/I$ satisfies the Serre condition $(R_1)$, hence $R/I$ is a normal domain.

\smallskip

(c) The estimates in (a) imply that $I$ satisfies the condition $(F_1)$ (or $G_{\infty}$).
Therefore, it is of linear type (see \cite{HSV1}).

\smallskip

(d) The assertion could possibly be derived from the methods of \cite{Nam}, but one can give a direct argument in the present
situation.
By part (c), $I$ is of linear type. Since $I$ is strongly Cohen--Macaulay (\cite[Theorem 2.1(a)]{AVRAMOVHERZOG}) then
the Rees algebra of $I$ is Cohen--Macaulay (\cite[Theorem 9.1]{HSV1}), and hence so is the associated graded ring of $I$.
On the other hand,  by part (b) the ideal $I$ is prime. By \cite[Proposition 3.2 (1)]{EH}, the assertion is equivalent to having
$$ \ell_P(I)\leq \max\{\hht P-1, \hht I\},$$
for every prime ideal $P\supset I$.
We may assume that $\hht P\geq 3$ since $I$ is a height $2$ prime. Therefore, we have to show that $\ell_P(I)\leq \hht P-1$.
If $P=(\XX)$ the result is clear since $\ell_{(\XX)}\leq \mu(I)=m\leq n-1=\hht (\XX)-1$.
Therefore, we may assume that $P\subsetneq (\XX)$, hence $\hht P\leq n-1$.
In particular, $I_1(\mathcal{C})\not\subset P$, so consider the index
$t_0:=\max\{1\leq s\leq m-2\,|\, I_s(\mathcal{C})\not\subset P\}$.
Pick a $t_0$-minor  $\Delta$ of $\mathcal{C}$ not contained in $P$, so that, in particular, $R_P$ is
a localization of the ring of fractions $R_{\Delta}=R[\Delta^{-1}]\subset k(\XX)$.
By a standard row-column elementary operation procedure, there is an $(m-t_0)\times(m-t_0-1)$ matrix $\widetilde{\mathcal{C}}$ over
$R_P$ such that $I_P=I_{m-1-t_0}(\widetilde{\mathcal{C}}).$

Now, we have $\hht I_{t_0+1}(\mathcal{C})-1\geq m-t_0$  from item (a) and, since $I_{t_0+1}(\mathcal{C})\subset P$ by definition of $t_0$,
it follows that $m-t_0\leq \hht P-1$.
Therefore
\begin{equation*}
\ell_P(I)=\ell (I_{m-1-t_0}(\widetilde{\mathcal{C}}))\leq\min\{\mu(I_{m-1-t_0}(\widetilde{\mathcal{C}})), \hht P\}
=\min\{m-t_0, \hht P\}\leq\hht P-1. \quad\quad\quad\quad\mbox{\rm \qed}
\end{equation*}

\subsection{Generic sub-Hankel matrices}

In this part we consider a degeneration of the generic Hankel (i.e., $1$-catalectic) matrix $\mathcal{H}$,
in which a lower corner of suitable size has its entries replaced  by zeros.
A version of this model for square matrices has been introduced in \cite{CRS}
(see also \cite{maron}) in connection to the construction of homaloidal determinants.

Let $R=k[X_1,\ldots,X_n]$ and let $m$ be such that $4\leq m+1\leq n\leq 2(m-1)$.
Set

\medskip

{\small
\begin{equation}\label{sub_hankel}
\mathcal{SH}=\left(
\begin{matrix}
X_{1}&\ldots &X_{n-m+1}&X_{n-m+2}&X_{n-m+3}& X_{n-m+4}&\ldots & X_{m-2}&X_{m-1}\\
X_{2}&\ldots &X_{n-m+2} & X_{n-m+3}&X_{n-m+4} & X_{n-m+5}&\ldots & X_{m-1}&X_{m}\\
\vdots&   &\vdots&\vdots & \vdots &  \vdots & &  \vdots & \vdots\\
X_{n-m+2}&\ldots &X_{2n-2m+2)}&X_{2n-2m+3}&X_{2n-2m+4}& X_{2n-2m+5}&\ldots &X_{n-1}&X_{n}\\[2pt]
\hline\\[-4pt]
X_{n-m+3}& \ldots &X_{2n-2m+3}&X_{2n-2m+4}&X_{2n-2m+5}&X_{2n-2m+6} &\ldots &X_n&0\\[4pt]
X_{n-m+4}& \ldots & X_{2n-2m+4}&X_{2n-2m+5}&X_{2n-2m+6}&X_{2n-2m+7} &\ldots & 0&0\\[6pt]
\vdots&   &\vdots&\vdots & \vdots &  \vdots & &  \vdots & \vdots\\[8pt]
X_{m-1}&\ldots &X_{n-1}& X_{n}& 0 &0 &\ldots & 0&0\\[4pt]
X_{m}&\ldots &X_{n}&0 &0& 0& \ldots &0&0
\end{matrix}
\right)
\end{equation}
}

Note that $n=2(m-1)$ is the case of the ordinary Hankel matrix.
This model has the following properties:

\begin{Theorem}\label{the_subs}
Let $m\geq 3$ and $n\leq 2(m-1)$.
Set $I:=I_{m-1}(\mathcal{SH})$.
Then:
\begin{enumerate}
\item[{\rm (a)}] $\hht(I_t(\mathcal{SH}))\geq m-t+2$ for $1\leq t\leq m-2$ and $I$ is a height $2$ prime ideal.
\item[{\rm (b)}] $R/I$ is normal if and only if $n\geq m+2$
\item[{\rm (c)}] $I$ is an ideal of linear type.
\item[{\rm (d)}] $I$ is normally torsionfree.
\end{enumerate}
\end{Theorem}
\demo
(a) The proof is similar to the argument used in the case of the catalecticant,
but there are differences due to the presence of zeros.
As before, the case $t=1$ is immediate.
Next we let $[\mathcal{SH}]_t$ denote the submatrix of $\mathcal{SH}$ with the first $t$ columns,
for values of $t$ in the range $2\leq t\leq m-2$.
For $t\leq n-m+1$ the matrix $[\mathcal{SH}]_t$ is an ordinary $m\times t$ Hankel matrix.
Therefore, it is $1$-generic, hence its
$t$-minors (maximal minors) generate a prime ideal of codimension $\geq m-t+1$ (also directly by the
observation after \cite[Proposition 4.3]{Eisenbud2}).
For $t > n-m+1$, one has
as follows
\begin{equation*}
k[X_1,\ldots, X_{n}]/I_t([\mathcal{SH}]_t)\simeq k[X_1,\ldots, X_n,X_{n+1},\ldots, X_{n+1+s}]/(X_{n+1},\ldots, X_{n+1+s},I_t([\mathcal{H}]_t)),
\end{equation*}
where $s=t-(n-m+1)$ and $[\mathcal{H}]_t$ is the Hankel matrix of size $m\times t$.
Since $s\leq m-2$ (because $t<n$) then we can use
\cite[Theorem 1]{Eisenbud1} to deduce again that $I_t([\mathcal{SH}]_t)$ is a prime ideal of codimension $m-t+1$.
The second assertion of (a) can be deduced similarly.

Now, to complete the argument it suffices to show that the inclusion $I_t([\mathcal{SH}]_t)\subset I_t(\mathcal{SH})$ is proper
for $2\leq t\leq m-2$, since the first of these ideals is prime.
This is clear if $t\leq 2(m-1)-n$ and $2(m-1)>n$ since a non-trivial power of $X_n$ is a $t$-minor in $I_t(\mathcal{SH})$
that cannot lie in the prime $I_t([\mathcal{SH}]_t)$ for $t\geq 2$.
Thus, assume that $t > 2(m-1)-n$. The argument is akin to the one used in the proof of the catalectic case.
 Namely, let $\Delta$ stand for the lower-rightmost $t$-minor of $\mathcal{SH}$
 {\small
 $$
 \left(\begin{array}{cccccccccccccccccc}
 X_{n-(t-r)-(t-1)}&X_{n-(t-r)-(t-1)+1}&\ldots&X_{n-(t-1)}&\ldots&X_{n-(t-r)-1}&X_{n-(t-r)}\\
 X_{n-(t-r)-(t-1)+1}&X_{n-(t-r)-(t-1)+2}&\ldots&X_{n-(t-1)+1}&\ldots&X_{2m-t}&X_{n-r}\\
 \vdots&\vdots&\ddots&\vdots&\ddots&\vdots\\
 X_{n-(t-1)}&X_{n-(t-1)+1}&\ldots&X_{n-(t-1)+(t-r)}&\ldots&X_{n-1}&X_n\\
 X_{n-(t-1)+1}&X_{n-(t-1)+2}&\ldots&X_{n-(t-1)+(t-r)+1}&\ldots&X_n&0\\
 \vdots&\vdots&\ddots&\vdots&\ddots&\vdots&\vdots\\
 X_{n-(t-r)}&X_{n-(t-r)+1}&\ldots&X_n&\ldots&0&0
 \end{array}\right)
 $$}
 where $r=2(m-1)-n$, which is the number of times $X_n$  appears on the original matrix.
By Laplace along the first row, it obtains

$$\Delta=X_{n-(t-r)-(t-1)}\left|\begin{array}{cccccccccccccccccc}
X_{n-(t-r)-(t-1)+2}&\ldots&X_{n-(t-1)+1}&\ldots&X_{2m-t}&X_{n-r}\\
\vdots&\vdots&\ddots&\vdots&\ddots&\vdots\\
X_{n-(t-1)+1}&\ldots&X_{n-(t-1)+(t-r)}&\ldots&X_{n-1}&X_n\\
X_{n-(t-1)+2}&\ldots&X_{n-(t-1)+(t-r)+1}&\ldots&X_n&0\\
\vdots&\ddots&\vdots&\ddots&\vdots&\vdots\\
X_{n-(t-r)+1}&\ldots&X_n&\ldots&0&0
\end{array}\right| + H,$$
where $H$ does not involve $X_{n-(t-r)-(t-1)}$ as the latter occurs only once on the matrix.

Now, by induction the determinant multiplying $X_{n-(t-r)-(t-1)}$ in this expression is of the form  $G\cdot X_n^{r}$, with $G$ a nonzero polynomial of degree  $(t-1)-r$.  This entails:

$$\Delta= (X_{n-(t-r)-(t-1)}G+F)X_n^r+(\mbox{terms of degree less than} \;r\; \mbox{in}\; X_n)$$
where $F$  (possibly vanishing) comes from terms in $H$, hence, in particular, does not involve  $X_{n-(t-r)-(t-1)}.$ From this, $X_{n-(t-r)-(t-1)}G+F\neq 0.$
On the other hand, any $t\times t$ submatrix of
 $\mathcal{SH}_t$ has at most $r-1$ entries equal to   $X_n$. 
 It follows that the $X_n$-degree on any $t$-minor generating  $I_t(\mathcal{SH}_t)$ is strictly less than $r.$ 
 
 Thus, we are led to conclude that $\Delta\notin I_t(\mathcal{SH}_t),$
 as claimed.

\smallskip

(b) We first show that $I$ does not satisfy $(R_1)$ if $n=m+1$ (the lowest possible value).
For this, consider the height $3$ prime  $P=(X_{n-2}, X_{n-1},X_n).$
Clearly, $I\subset P$ by direct inspection on the shape of the matrix.
Note that the upper left $(n-3)$-minor of $\mathcal{SH}$ has the form $X_{n-3}^{n-3}+q$ where $q\in P$, hence does not belong to $I$.
After appropriate row/column operations,  we se that $I_P=(\Delta_{n-2},\Delta_{n-1})$, where $\Delta_i$ denotes $(n-2)$-minor
 of $\mathcal{SH}$ obtained by omitting the $i$th row. We claim that  $R_P/I_P$ is not regular.
 For this, it suffices to show that $\Delta_{n-2}\in P^{2}.$
But

{\scriptsize \begin{eqnarray*}
\Delta_{n-2}=(-1)^{n-1}X_{n-1}\det\left(\begin{array}{cccccccc}X_2&X_3&\ldots &X_{n-2}\\
X_3&X_4&\ldots&X_{n-1}\\ \vdots&\vdots&\vdots&\vdots\\ X_{n-2}&X_{n-1}&\ldots&0\end{array}\right)
+ (-1)^{n}X_{n}\det\left(\begin{array}{cccccccc}X_1&X_3&\ldots &X_{n-2}\\ X_2&X_4&\ldots&X_{n-1}\\
 \vdots&\vdots&\vdots&\vdots\\ X_{ n-3}&X_{n-1}&\ldots& 0\end{array}\right)
\end{eqnarray*}}

Note that the two determinants have same last column and the nonzero entries on the column are $X_{n-2},X_{n-1},X_n$.
Therefore, expanding these determinants along their last column  clearly shows the claim.

\smallskip

Conversely, suppose now that $n\geq m+2$.
We will use use a result of \cite{Eisenbud2}.
For this, first switch to the transpose $M$ of $\mathcal{SH}$ in order to conform with the notation in
\cite{Eisenbud2}.
Next let $M'$ denote the transpose of the Hankel matrix of size $(m-1)\times m$.
Thinking of $M\subset M'$ as the respective $k$-subspaces spanned by the entries, one has
${\rm codim}_{M'}M=2(m-1)-n$.
Our hypothesis implies that ${\rm codim}_{M'}M\leq m-4$.
Now, $M'$ is $1$-generic, hence the singular locus of ${\rm Proj}(R/I)$ is
contained in the union of ${\rm Proj}(R/I_{m-2}(\mathcal{SH}))$ and a certain set of codimension
at least $m-2-{\rm codim}_{M'}M$ in ${\rm Proj}(R/I)$ -- 
according to the discussion immediately following the statement of \cite[Theorem 2.1 (3)]{Eisenbud2} and 
its proof in \cite[Proposition 3.1 and The completion of the proof of Theorem 2.1]{Eisenbud2}.
The first of these two has codimension $\geq 2$ in ${\rm Proj}(R/I)$ by the estimates of item (a).
As for the second, its codimension in ${\rm Proj}(R/I)$ is now at least $m-2-{\rm codim}_{M'}M\geq m-2-(m-4)=2$.
This shows that $R/I$ satisfies Serre's property $(R_1)$, hence is normal.

\smallskip

(c) and (d) are proved exactly the same way as in Theorem~\ref{Perfect_catalecticant}
\qed

\subsection{Semi-Hankel matrices}

Let again $n\leq 2(m-1)$, this time around with $3\leq m\leq n$, thus allowing for the equality $n=m$
(while $n\geq m+1$ was stipulated in the sub-Hankel case.)
The model of this part is in a sense an ancestral of the previous sub-Hankel model, in which one specializes certain entries
of the $m\times (m-1)$ generic Hankel  matrix to a few independent linear forms:

{\small
\begin{equation*}
\widetilde{\mathcal{H}}_{n,m}\left(
\begin{matrix}
X_{1}&\ldots &X_{n-m+1}&X_{n-m+2}&X_{n-m+3}&\ldots & X_{m-2}&X_{m-1}\\
X_{2}&\ldots &X_{n-m+2} & X_{n-m+3}&X_{n-m+4} &\ldots & X_{m-1}&X_{m}\\
\vdots&   &\vdots&\vdots & \vdots &   &  \vdots & \vdots\\
X_{n-m+2}&\ldots &X_{2n-2m+2)}&X_{2n-2m+3}&X_{2n-2m+4}&\ldots &X_{n-1}&X_{n}\\[2pt]
\hline\\[-4pt]
X_{n-m+3}& \ldots &X_{2n-2m+3}&X_{2n-2m+4}&X_{2n-2m+5} &\ldots &X_n&\ell_1\\[4pt]
X_{n-m+4}& \ldots & X_{2n-2m+4}&X_{2n-2m+5}&X_{2n-2m+6} &\ldots & \ell_1&\ell_2\\[6pt]
\vdots&   &\vdots&\vdots & \vdots  & &  \vdots & \vdots\\[8pt]
X_{m-1}&\ldots &X_{n-1}& X_{n}& \ell_{1}  &\ldots & \ell_{2(m-1)-n-2}&\ell_{2(m-1)-n-1}\\[4pt]
X_{m}&\ldots &X_{n}&\ell_{1} &\ell_{2}& \ldots &\ell_{2(m-1)-n-1}&\ell_{2(m-1)-n}
\end{matrix}
\right)
\end{equation*}
}

\noindent where $\{\ell_{1}, \ell_{2}, \ldots ,\ell_{2(m-1)-n-1}, \ell_{2(m-1)-n}\}$ are independent 
linear forms in $k[X_1,\ldots, X_n]$, an assumption that makes sense since $2(m-1)-n\leq n$.

\begin{Theorem}\label{main_semi_hankel}
Let $I:=I_{n-1}(\widetilde{\mathcal{H}}_{n,m})$. Then:
\begin{enumerate}
\item[{\rm (a)}] $\hht(I_t(\widetilde{\mathcal{H}}_{n,m}))\geq m-t+2$ in the range $1\leq t\leq m-2$, while  $\hht(I)=2$
\item[{\rm (b)}] $R/I$ is a normal domain.
\item[{\rm (c)}] $I$ is an ideal of linear type.
\item[{\rm (d)}]  {\rm (char$(k)=0$)} For every $r\geq 0$ such that $I^{(r)}\neq I^{r}$,  the $R$-module $I^{(r)}/I^{r}$ is
$(\XX)$-primary.
\end{enumerate}
\end{Theorem}
\demo
(a) We use same strategy as before, by considering the submatrix $[\widetilde{\mathcal{H}}_{n,m}]_t$ of $\widetilde{\mathcal{H}}_{n,m}$
formed with the first $t\leq m-1$ columns.
Since this matrix specializes from the full generic $m\times t$ Hankel matrix -- which is $1$-generic -- modulo a regular sequence of $1$-forms
of cardinality  $t-(n-m+1)\leq t-2$, we can apply
\cite[Theorem 2.1 and Corollary 3.3]{Eisenbud2}, by which one has the following properties:
\begin{enumerate}
\item[(i)] $\hht I_t([\widetilde{\mathcal{H}}_{n,m}]_t)=m-t+1,$ for $2\leq t\leq m-1$.
\item[(ii)] $I_t([\widetilde{\mathcal{H}}_{n,m}]_t)$ is a prime ideal.
\item[(iii)] The ideal $I_{t-1}([\widetilde{\mathcal{H}}_{n,m}]_t)/I_t([\widetilde{\mathcal{H}}_{n,m}]_t)$ defines the
singular locus of $R[\ZZ]/I_t([\widetilde{\mathcal{H}}_{n,m}]_t)$
\end{enumerate}
Consider the $t$-minor $\Delta$ of the lower rightmost corner of $\widetilde{\mathcal{H}}_{n,m}$.
The argument is the same as of the proof of Theorem~\ref{the_subs} (a): we may assume that $t\geq 2(m-1)-n$,
hence $\Delta$ has a leading term in $X_n$ coming from the anti-diagonal with $X_n$'s throughout -- note that the terms involving
the linear forms can only have a smaller degree in $X_n$.
Again, we are led to conclude by the same token as before that $\Delta\notin I_t([\widetilde{\mathcal{H}}_{n,m}]_t)$ for $2\leq t\leq m-2$.

This proves that $\hht(I_t(\widetilde{\mathcal{H}}_{n,m}))\geq m-t+2$ in the range $1\leq t\leq m-2$, and also that  $\hht(I)=2$.

\smallskip

(b) Using the result of (iii) in the case of $t=m-1$, we know that the singular locus of $R/I$ is defined by the ideal
$I_{m-2}(\widetilde{\mathcal{H}}_{n,m})/I_{m-1}(\widetilde{\mathcal{H}}_{n,m})$.
By (a) the latter has codimension $\geq m-(m-2)+2-2=2$, hence $R/I$ has the property $(R_1)$.
Therefore, $R/I$ is a (Cohen--Macaulay) normal domain.

\smallskip

(c) As previously remarked, (a) implies condition $(F_1)$.
Since $R/I$ is Cohen--Macaulay of codimension $2$, $I$ is of linear type.

\smallskip

(d)
Fixing an $r\geq 0$, suppose that $I^{(r)}/I^r\neq \{0\}$.
The assertion is equivalent to saying that a power of $(\XX)$ annihilates $I^{(r)}/I^r$
i.e., that ${I^{(r)}}_P={I^r}_P$ for every prime $P\neq (\XX)$.
Letting $r\geq 0$ run, this is in turn equivalent to claiming that the associated graded ring gr$_I(R)$ is torsionfree
over $R/I$ locally on the punctured spectrum Spec$(R)\setminus (\XX)$.

Thus, let $P\neq (\XX)$ be a prime containing $I$.
By (a), as already pointed out, $I$ and hence, also, $I_P$ satisfies the condition $(F_1)$.
As in the proof of Theorem~\ref{Perfect_catalecticant} (d), we know that the associated graded ring gr$_{I_P}(R_P)$
is Cohen--Macaulay. Therefore, by the same token and since $\hht I=2$, one has to show the local estimates
$$\ell_Q(I)=\ell_{Q_P}(I_P)\leq \hht (Q_P)-1= \hht Q-1,$$
for every prime $Q\subset P$.

Fixing such a prime $Q$, set $t_0:=\max\{1\leq s\leq m-2\,|\, I_s(\widetilde{\mathcal{H}}_{n,m})\not\subset Q\}$ - again $t_0$ makes sense
since $I_1(\widetilde{\mathcal{H}}_{n,m})\not\subset Q$.
Trading $Q$ for $P$ in the proof of Theorem~\ref{Perfect_catalecticant} (d), the rest of the argument is literally the same.
\qed

\medskip

In the case where $n=m$, we gather some geometric information:

\begin{Proposition}
Let $I:=I_{n-1}(\widetilde{\mathcal{H}}_{n,n})$. Then:
\begin{enumerate}
\item[{\rm (i)}] $I^{(\ell)}=I^{\ell}$ for $1\leq \ell\leq n-2$
\item[{\rm (ii)}] The rational map $\mathfrak{G}:\pp^{n-1}\dasharrow\pp^{n-1}$ defined by the $(n-1)$-minors of
$\widetilde{\mathcal{H}}_{n,n}$ is a Cremona map.
\item[{\rm (iii)}] {\rm (char$(k)=0$)} The symbolic Rees algebra
$\mathcal{R}^{(I)}$ of $I$ is a Gorenstein normal domain such that $\mathcal{R}^{(I)}=R[It,Dt^{n-1}]$,
where $D$ is the source inversion factor of the Cremona map defined by the $(n-1)$-minors of
$\widetilde{\mathcal{H}}_{n,n}${\rm ;}  moreover, $D$ coincides with the Jacobian determinant of the
 $(n-1)$-minors of $\widetilde{\mathcal{H}}_{n,n}$.
 \end{enumerate}
\end{Proposition}
\demo
(i)  By Theorem~\ref{main_semi_hankel} (c),  $I$ is an ideal of linear type; in particular,
it satisfies property $(G_{\infty})$, i.e., $\mu(I_p)\leq \hht P$ for every prime $P\subset R$.
According to \cite[Theorem 5.1]{TCHERNEV}, under this condition, for every $1\leq\ell\leq n-1$ the complex
$$\mathcal{K}_{\ell}\;:\;0\rar F_{\ell}\rar F_{\ell-1}\rar\ldots\rar F_1\rar F_0\rar0$$
is a free resolution of $I^{\ell}$, where $F_i:=\bigwedge^{i}R^{n-1}\otimes_R \mathcal{S}_{(n-1)-i}(R^{n})$
and  $d:F_i\rar F_{i-1}$ is defined through
$$d(e_1\wedge\ldots\wedge e_i\otimes g):=\sum_{l=1}^{i}e_1\wedge\ldots\wedge\widehat{e_l}
\wedge\ldots\wedge e_i\otimes\phi(e_l)g.$$
It follows that, in the range $1\leq\ell\leq n-2$,  the $R$-module $(R/I^{\ell})$ has homological dimension
 at most $n-1$, hence $(\XX)$ is not an associated prime thereof.
By induction using $\ass(R/I^{\ell})\subset \ass(I^{(\ell)}/I^{\ell})\cup\ass(R/I^{(\ell)})$ and drawing upon
Theorem~\ref{main_semi_hankel}, we are done.

\smallskip

(ii) This follows from Theorem~\ref{main_semi_hankel} (a) and (c) via \cite[Proposition 3.4]{AHA}.

\smallskip

(iii)
The symbolic Rees algebra
$\mathcal{R}^{(I)}$ of $I$ is a Gorenstein ring; indeed, it is a quasi-Gorenstein Krull
domain since $I$ is a codimension $2$ prime ideal (\cite{ST}).
On the other hand, by the proof of \cite[Corollary 3.4 (b)]{dual}, $\mathcal{R}^{(I)}$ is
finitely generated since one has an isomorphism
$\mathcal{R}^{(I)}\simeq \mathcal{R}(I)[t^{-1}]=R[It,t^{-1}]$. Moreover, the latter is Cohen--Macaulay since
$\mathcal{R}(I)$ is Cohen--Macaulay -- same argument as in the proof of \cite[Proposition 2.9 (b)]{zaron}.
It follows that $\mathcal{R}^{(I)}$ is a Gorenstein normal domain and that it is generated over
the Rees algebra $R[It]$ by one single element.
It remains to identify this element in the form $Dt^d$, for some $D\in R$ and a uniquely determinde exponent $d$.

To go about this, we use (ii), namely, let $\mathfrak{d}_1,\ldots,\mathfrak{d}_n\in k[\YY]$ be forms
of the same degree, with $\gcd =1$, defining the inverse map and let $D\in R$ denote the corresponding
source inversion factor.
Write $J=(\mathfrak{d}_1,\ldots,\mathfrak{d}_n)\subset k[\YY]$.
By definition, one has
$$D=\mathfrak{d}_i(\Delta_{1},\ldots,\Delta_{n})/X_i,\, 1\leq i\leq n,$$
where $\mathbf{\Delta}:=
\{\Delta_{1},\ldots,\Delta_{n}\}$ are the (signed) minors generating $I$.
But under a Cremona map, the two Rees algebras $\mathcal{R}_R(I)=R[It]\subset R[t]$ and
$\mathcal{R}_{k[\YY]}(J)=k[\YY][Ju]\subset k[\YY][u]$ get identified by a $k$-isomorphism
that maps $Y_i\mapsto \Delta_i t$ and $X_i\mapsto \mathfrak{d}_iu$ (see. e.g., \cite[Proposition 2.1]{bir2003}).
Then $D$ is identified with
$\mathfrak{d}_1/X_1$ in the common field of fractions.
Using (i) above,  the symbolic algebra is generated
by $It$ and $Dt^{n-1}$ as a consequence of
\cite[Corollary 7.4.3 (b)]{Wolmbook}.

The additional statement reads out of Proposition~\ref{inversion_factor_is_determinant}.
\qed


\noindent {\bf Authors' addresses:}

\medskip

\noindent {\sc Zaqueu Ramos}, Departamento de Matem\'atica, CCET, Universidade Federal de Sergipe\\
49100-000 S\~ao Cristov\~ao, Sergipe, Brazil\\
{\em e-mail}: zaqueu.ramos@gmail.com\\

\noindent {\sc Aron Simis},  Departamento de Matem\'atica, CCEN, Universidade Federal
de Pernambuco\\
 50740-560 Recife, PE, Brazil.\\
{\em e-mail}:  aron@dmat.ufpe.br

\end{document}